\documentclass[10pt]{article}
\usepackage{amssymb,amsmath, amsthm, amsfonts}
\usepackage{mathrsfs,color}
\usepackage[all]{xypic}
\usepackage{graphicx}
\usepackage{tikz,color} 
\usetikzlibrary{calc,matrix,arrows,chains,positioning,scopes,decorations.pathmorphing,backgrounds,fit}
\newcommand{\luk}{\L u\-ka\-si\-e\-w\-icz}
\usepackage{yfonts}

\usepackage{enumitem}
\usepackage{adjustbox}
\usepackage{thmtools,mathtools}

\usepackage{nameref,
cleveref}

\declaretheorem[numberwithin=section,style=theorem,refname={theorem,theorems},  Refname={Theorem,Theorems}, name=Theorem]{theorem}
\declaretheorem[sibling=theorem,style=remark,refname={remark,remarks}, Refname={Remark,Remarks},name=Remark]{remark}
\declaretheorem[sibling=theorem,style=theorem,refname={proposition,proposition},  Refname={Proposition,Propositions},name=Proposition]{proposition}
\declaretheorem[sibling=theorem,style=theorem, name=Lemma, refname={lemma,lemmata},  Refname={Lemma,Lemmata}]{lemma}
\declaretheorem[sibling=theorem,style=theorem,refname={corollary,corollaries},  Refname={Corollary,Corollaries},name=Corollary]{corollary}

\declaretheorem[numbered=no,style=definition,refname={notation,notations},  Refname={Notation,Notations},name=Notation]{notation}

\declaretheorem[sibling=theorem,style=definition,refname={definition,definitions},  Refname={Definition,Definitions},name=Definition]{definition}

\newcommand{\bdfn}{\begin{definition}}
\newcommand{\edfn}{\end{definition}}
\newcommand{\bthm}{\begin{theorem}}
\newcommand{\ethm}{\end{theorem}}
\newcommand{\bprop}{\begin{proposition}}
\newcommand{\eprop}{\end{proposition}}
\newcommand{\bcor}{\begin{corollary}}
\newcommand{\ecor}{\end{corollary}}
\newcommand{\blem}{\begin{lemma}}
\newcommand{\elem}{\end{lemma}}
\newcommand{\bfact}{\begin{remark}}
\newcommand{\efact}{\end{remark}}
\newcommand{\bex}{\begin{example}}
\newcommand{\eex}{\end{example}}
\newcommand{\bnot}{\begin{notation}}
\newcommand{\enot}{\end{notation}}
\newcommand{\ten}{\otimes}
\newcommand{\quot}[2]{{\raisebox{.2em}{$#1$}\left/\raisebox{-.2em}{$#2$}\right.}}

\tikzset{node distance=2.3cm, auto}

\begin{document}

\title{A general view of the algebraic semantics of  \L ukasiewicz logic with product}
\author{Serafina Lapenta\footnote{Corresponding author}\\ 
{\small Department of Mathematics, University of Salerno,}\\ 
{\small  Via Giovanni Paolo II, 132 Fisciano (SA), Italy}\\
 {\small slapenta@unisa.it}
 \and
Ioana Leu\c stean \\
{\small Department of Computer Science,} \\
{\small Faculty of Mathematics and Computer Science, University of Bucharest,}\\
{\small Academiei nr.14, sector 1, C.P. 010014,  Bucharest, Romania}\\   
{\small ioana.leustean@unibuc.ro}
}
\date{}
\maketitle

\begin{abstract}
This paper aims at connecting the various classes that provide an algebraic semantics for three different conservative expansions of \L ukasiewicz logic, using algebraic and categorical theoretic techniques. We connect such classes of algebras by adjunctions, using the tensor product of MV-algebras and defining the \textit{tensor PMV-algebra} of a semisimple MV-algebra, inspired by the construction of the tensor algebra of a vector space. We further apply the main results to prove amalgamation properties and, via categorical equivalence, we transfer all results to the framework of lattice-ordered groups.\\
\textit{Keywords:} 
MV-algebra, tensor product, scalar extension property, tensor algebra, amalgamation property.
\end{abstract}

\section*{Introduction}
In the late Fifties C.C. Chang introduced the notion of MV-algebra in the pursuit of a simpler proof of completeness for the infinite-valued \luk\ logic, inspired by the theory of lattice-ordered groups. Almost thirty years later, D. Mundici proved that MV-algebras are categorically equivalent to Abelian lattice-ordered groups with a strong unit. This valuable result provided the stepping stone needed for an investigation of MV-algebras in the spirit of the theory of groups, rings and  algebras.

Indeed, stemming out from such a relation with groups, several expansions of MV-algebras have been defined: one can find the notion of PMV-algebras (MV-algebras endowed with a ``ring-like" product), Riesz MV-algebras (MV-algebras endowed with a scalar multiplication), \textit{f}MV-algebras (MV-algebras endowed with both ring-like product and scalar product), MV-modules (MV-algebras endowed with the action of a PMV-algebra on them). The initial motivation came from logic: since the standard MV-algebra on $[0,1]$ is closed to the real product, it was natural to look for complete theories of $[0,1]$ endowed with \L ukasiewicz operations and powerful enough to axiomatize the product of real numbers. 

All of the above-mentioned structures have been investigated from a category-theoretical and a universal-algebraic point of view, obtaining important results: categorical equivalences with rings, vector lattices, \textit{f}-algebras and lattice-ordered modules; geometrical dualities with suitable categories of polyhedra \cite{MunBook,DMVNotes, RMV-DiNola, Lpoly}; functional representation of free algebras and normal forms theorems \cite{LLNF}; a suitable probability theory through the notion of states \cite{MunBook,FlaKroupa-HB}.  

PMV-algebras, Riesz MV-algebras and \textit{f}MV-algebras are, of course, deeply related to each-other. Trivially, all of them have an MV-algebra reduct,  moreover any \textit{f}MV-algebra has a PMV-algebra reduct as well as a Riesz MV-algebra reduct. Thus, it is natural to ask if these forgetful functors can be reversed allowing on the one end to gain a deeper understanding of all classes of algebras involved and on the other end to obtain new knowledge on the logical systems attached to each class of algebras via algebraic and categorical techniques.

With this two main objectives in mind, we started this investigation in \cite{LLTP1} where, using the tensor product of semisimple MV-algebras \cite{Mun}, we connected semisimple MV-algebras with semisimple Riesz MV-algebra and semisimple PMV-algebras with semisimple \textit{f}MV-algebras by categorical adjunctions. We provided these connections by displaying two pairs of adjunctions and by proving the \emph{scalar extension properties} for semisimple algebras, properties that were fundamentally used as main algebraic tool.

In this paper we conclude the work of \cite{LLTP1}, giving a clear and comprehensive picture of the relations between the aforementioned structures through the classical construction of the \textit{tensor algebra} of a vector space. We define the \textit{tensor PMV-algebra} of a semisimple MV-algebra and apply the same construction to a Riesz MV-algebra in order to obtain an \textit{f}MV-algebra. In Section \ref{sec:CatSett} we lift these results to a categorical level, obtaining two pairs of adjoint functors, connecting semisimple MV-algebras with semisimple PMV-algebras and semisimple Riesz MV-algebras with semisimple \textit{f}MV-algebras. Figure 8 frames our results in the existing literature and provides the complete graph of relationships between the different classes of semisimple MV-algebras with products.

In Section \ref{amalg:PMV} and Section \ref{sec:free} we prove how the results from the previous sections can have an impact of the logical systems attached to the considered classes of algebras, giving at the same time two examples of the value of our constructions by itself. In Section \ref{amalg:PMV} we  prove that semisimple PMV-algebras, semisimple Riesz MV-algebras and semisimple \textit{f}MV-algebras have the  amalgamation property, while in Section \ref{sec:free} we obtain three characterizations of the free objects in the categories involved. These results have a logical taste, indeed we remark that amalgamation in varieties is closely related to  interpolation in logic and  free algebras in this framework are isomorphic to the Lindenbaum-Tarski algebras of the corresponding logical system.

Finally, in Section \ref{sec:lgroups} we transfer all results to lattice-ordered structures via categorical equivalence. In this way we obtain new properties for groups, rings, vector lattices and algebras directly from MV-algebras.

We remark that all results are presented for semisimple MV-algebras and Archimedean lattice-ordered groups, choice which is thoroughly motivated in Section \ref{sec:SEP} and it is more than adequate when one deals with algebraic logic.

\section{Preliminaries}
\subsection{An overview of MV-algebras with product and \L ukasiewicz logic}\label{sec:intrMV}
In this section we provide a short overview of MV-algebras and their expansions. Nonetheless, we urge the interested reader to consult \cite{CDM,MunBook,DiNLeu-HB} for an in-depth treatment.

MV-algebras are the algebraic counterpart of the infinite-valued \luk\ logic. They are structures $(A, \oplus, \ ^*, 0)$, where $(A, \oplus ,0)$ is a commutative monoid, $\ ^*$ is an involution and the equations $(x^*\oplus y)^*\oplus y=( y^*\oplus x)^*\oplus x$  and $x \oplus 0^*=0^*$ are satisfied for any $x,y,z\in A$. Any MV-algebra can be endowed with a lattice order and  the standard model is the unit interval $[0,1]$ with $x\oplus y= \min (x+y,1)$ and $x^*=1-x$. The variety of MV-algebras is generated by the standard model and since $[0,1]$ is closed with respect to the real product, a fruitful  research direction  is the study of MV-algebras enriched with a product operation \cite{DiND,Mon+,MoPa,HC,LeuRMV}, which can be a  binary  operation, a scalar multiplication or a combination of both. 

If $A,B, C$ are MV-algebras, a function $\omega:A\to B$  is called {\em linear}  if $\omega(x\oplus y)=\omega(x)\oplus\omega(y)$ whenever $x\leq y^*$. A function $\beta:A\times B\to C$ is {\em bilinear} if  $\beta(-,y)$ and $\beta(x,-)$ are linear for any $x\in A$ and $y\in B$.  These notions allow us to define the whole algebraic hierarchy of MV-algebras with product in a uniform way, as described in Table 1\footnote{A more general definition of PMV-algebras and \textit{f}MV-algebras can be found in \cite{DiND,LLfMV} respectively. Nonetheless, the current paper focuses on unital and commutative structures.}.
\begin{center}
\normalsize{
\begin{tabular}{|c|l|}\hline
{\bf Structure} & {\bf Definition} \\ \hline\hline
 &   $(P,\oplus, ^*,0 )$ MV-algebra \\
   $(P,\oplus,\cdot, ^*, 0)$ &
 $\cdot : P\times P\to P$ bilinear,\\
unital and commutative  &  $x \cdot (y\cdot z)= (x \cdot y)\cdot z$,
\\
PMV algebra  \cite{DiND,Mon+} & $x\cdot 1=1\cdot x=x$\\ \hline
                                            &   $(R,\oplus, ^*,0 )$ MV-algebra \\
   $(R,\oplus, ^*, \{\alpha\mid\alpha\in [0,1]\}, 0)$ &
 $(\alpha, x)\mapsto\alpha x$ bilinear,\\
Riesz MV- algebra \cite{LeuRMV} &   $(\alpha \cdot \beta)x= \alpha (\beta x)$,
\\
& $1x=x$ \\ \hline
   $(A,\oplus,\cdot,  ^*, \{\alpha\mid\alpha\in [0,1]\}, 0)$& $(A,\oplus, \cdot, ^*,0 )$ unital PMV-algebra \\
 unital and commutative  &  $(A,\oplus,^*, \{\alpha\mid\alpha\in [0,1]\}, 0)$ R. MV-algebra  \\
 \textit{f}MV- algebra \cite{LLfMV}&  $\alpha (x\cdot y)=(\alpha x)\cdot y=x\cdot (\alpha y)$
 \\ \hline

\end{tabular}
\vspace*{0.2cm}

{Table 1. Algebraic hierarchy.}
}
\end{center}

The present investigation is centred on  the class of {\em semisimple} MV-algebras. Such algebras enjoy a crucial functional representation. Indeed, any semisimple MV-algebra is isomorphic to a separating MV-subalgebra of $[0,1]$-valued continuous functions defined over some compact Hausdorff space \cite{CDM}, namely the space of maximal ideals of the algebra. A PMV-algebra (Riesz MV-algebra or \textit{f}MV-algebra) is semisimple if  its MV-algebra reduct is semisimple\footnote{We remark that this holds in the case of unital PMV-algebras and unital \textit{f}-MV-algebras. More detail on semisimplicity for non-unital \textit{f}MV-algebras can be found in \cite{LLfMV}.}. 

For all structures defined in Table 1 it is possible to give an equational characterization. In other words, PMV-algebras, Riesz MV-algebras and \textit{f}MV-algebras are varieties, that we shall denote with $\mathbb{PMV}$, $\mathbb{RMV}$ and $\mathbb{FMV}$ respectively. As for MV-algebras, $\mathbb{RMV}$ is generated by its standard model $[0,1]$, where the scalar operation coincides with the product of real numbers, while in the case of PMV-algebras  the standard model generates a proper sub-class, more specifically, the class of \textit{semiprime} algebras, i.e. algebras defined by the quasi-identity ``$x^2=0$ implies $x=0$'' \cite{HC,Mon+} and such quasi-variety is denoted by $PMV^+$. A similar results is obtained in \cite{LLfMV}, where semiprime fMV-algebras are introduced and their quasi-variety is denoted by $FR^+$, while $FR$-algebras are elements of the variety generated by $[0,1]$.\footnote{Note that the result presented in \cite{LLfMV} was based upon a result in Birkhoff's Lattice Theory \cite{Birk}, for which we recently found a counterexample, see \cite[Example 11.54]{ce}. Therefore, despite the claim in \cite{LLfMV}, it is still not clear whether the quasi-variety $FR^+$ is generated by the standard \textit{f}MV-algebra $[0,1]$.}.\\

The natural hierarchy of lattice-ordered structures we have introduced in Table 1 is naturally reflected by appropriate logical systems. Indeed, in \cite{HC} the logics $P$\L \ and $P$\L$^\prime$ are defined and they have PMV-algebras and PMV$^+$-algebras respectively as models. In \cite{LeuRMV} the logical system $\mathcal{R}$\L \ is defined, and its models are Riesz MV-algebras. Both $P$\L \ and  $\mathcal{R}$\L\ are conservative extensions of \L ukasiewicz logic, and $P$\L$^\prime$ is obtained by $P$\L \ adding an appropriate deduction rule. Finally, in \cite{LLfMV} one can find the logical systems $\mathcal{FMVL}$ and $\mathcal{FMVL}^+$, whose models are $\textit{f}$MV-algebras and FR$^+$-algebras. We remark that the former is obtained extending the union of $P$\L\ and $\mathcal{R}$\L , while the latter extends the union of $P$\L$^\prime$ and $\mathcal{R}$\L .
%
%

By general results in Universal Algebra, free MV(Riesz, PMV$^+$, FR)-algebra $k$-generated exists and it is the subalgebra of  $([0,1])^{[0,1]^k}$ generated by the projection maps \cite{Gr}. Moreover, in each case, such free algebras  are the Lindenbaum-Tarski algebras (i.e. equivalence classes of formulas) of the corresponding logics. 

More specifically, McNaughton's theorem \cite{McN} states that the free $k$-generated MV-algebra is (up to isomorphism) the algebra of continuous 
functions from $[0,1]^k$ to $[0,1]$ that are piecewise
linear with integer coefficients. In other words, any element of the
Lindenbaum-Tarski algebra is a continuous function and there are
finitely many affine linear functions such that in any point of the
domain, it coincides with one of them. A
similar result holds for Riesz MV-algebras \cite[Theorem
11]{LeuRMV}. In this case the free object is the algebra of functions
which are piecewise linear with \emph{real} coefficients. In what follows we
will denote the free $k$-generated MV-algebra by $Free (k)$, while
$Free_R (k)$ denotes the free $k$-generated Riesz MV-algebra.  In the case of PMV$^+$-algebras and FR-algebras the characterization of the free object in terms of piecewise functions is an open problem and it is related to the long-standing Pierce--Birkhoff conjecture \cite{LLfMV}.\\

Finally, we recall that the category of MV-algebras is equivalent with the category of  Abelian lattice-ordered groups with strong unit \cite{Mun1}, {\em $\ell u$-groups} for short. If $(G, u)$ is an $\ell u$-group, the interval $[0,u]_G=\{ x\in G \mid 0\le x\le u\}$ (called \textit{unit interval}), is an MV-algebra with $ x\oplus y=u\wedge (x+y)$, $x^*=u-x$.  If $\mathbf{MV}$ is the category of MV-algebras and $\mathbf{auG}$ is the category of $\ell u$-groups equipped with morphisms that preserve the strong unit,  one can define a functor $\Gamma: \mathbf{auG}\rightarrow \mathbf{MV}$ by
$\Gamma (G,u)=[0,u]_G $ and $\Gamma(h)= h|_{[0,u_1]_{G_1}}$,
where  $(G,u)$ is an $\ell u$-group   and   $h: G_1\to G_2$ is a morphism in $\mathbf{auG}$ between the two $\ell u$-groups $(G_1,u_1)$ and $(G_2, u_2)$.  The functor $\Gamma$ establishes a categorical equivalence between  $\mathbf{auG}$ and $\mathbf{MV}$ \cite{Mun1}. 
Moreover, through $\Gamma$, semisimple MV-algebras correspond to Archimedean $\ell u$-groups. We shall denote by $\mathbf{MV_{ss}}$ the  full subcategory of semisimple MV-algebras and by $\mathbf{auG_a}$ the full subcategory of Archimedean $\ell u$-groups.

Extending $\Gamma$, similar equivalences are proved for:  PMV-algebras and a subclass of  lattice-ordered rings with strong unit (unital $\ell u$-rings, shortly);   Riesz MV-algebras and  Riesz spaces (vector lattices)  with strong unit; \textit{f}MV-algebras and \textit{f}-algebras with strong unit. The functors that give the equivalences are denoted by $\Gamma_{(\cdot)}$, $\Gamma_{\mathbb R}$ and $\Gamma_f$ respectively. See \cite{Birk,BP} for details on the above mentioned structures and \cite{DiND,LeuRMV,LLfMV} for details on the categorical equivalences. 

In Table 2, we set notations of all categories of semisimple and Archimedean structures involved in this investigation. We remark that unital and semisimple PMV-algebras and \textit{f}MV-algebras are commutative.

\begin{center}
\normalsize{
\begin{tabular}{|c|l|c|l}\hline
{\bf Category} & {\bf Objects} \\ \hline\hline

$\mathbf{uPMV_{ss}}$ & unital and semisimple PMV-algebras,\\
 $\mathbf{uR_{a}}$ & unital and Archimedean $\ell$-rings with strong unit, \\ \hline
$\mathbf{RMV_{ss}}$ & semisimple Riesz MV-algebras,\\
$\mathbf{uRS_a}$ & Archimedean Riesz spaces with strong unit,\\ \hline
$\mathbf{ufMV_{ss}}$ & unital and semisimple \textit{f}MV-algebras,\\
 $\mathbf{fuAlg_{a}}$ & unital and Archimedean \textit{f}-algebras with strong unit.\\ \hline
\end{tabular}
\vspace*{0.1cm}

{Table 2. Categories of MV-algebras and related $\ell$-structures.}

}
\end{center}
 
There are obvious forgetful functors between the above-defined categories and they commute with the $\Gamma$-type functors. 
A natural problem is to define appropriate left adjoints  for the forgetful functors. We started this investigation in \cite{LLTP1}, where the key tool was the {\em semisimple tensor product of MV-algebras} \cite{Mun} and, in particular, its {\em scalar extension property} \cite{LLTP1}, which will be discussed in the next subsection.

\subsection{The tensor product of lattice-ordered structures}\label{sec:intrTP}

The classical construction of a tensor product has been defined in the setting of lattice-ordered structures by several authors. In \cite{Mart}, the author defines a $\ell u$-bilinear function as a map
$\gamma: G\times H\rightarrow L$ between $\ell u$-groups $(G, u_G)$, $(H, u_H)$ and $(L, u_L)$  
such that $\gamma (x, -)$ and $\gamma (- , y)$ are homomorphisms of $\ell$-groups when $x$ and $y$ are positive and $\gamma (u_G, u_H) \le u_L$. 

Then, the tensor product is an $\ell u$-group $(G\ten_{\ell} H, u_G \ten_{\ell}u_H)$ together with an $\ell u$-bilinear map $ \gamma_{G,H}: G\times H \rightarrow G\ten_{\ell} H$  uniquely characterized, up to isomorphism, by a universal property with respect to $\ell$-groups \cite[Theorem 3.1]{Mart}. The map 
$ \gamma_{G,H}: G\times H \rightarrow G\ten_{\ell} H$ is defined by $\gamma_{G,H}(x,y)=x\ten_{\ell} y$.
The tensor product of Archimedean $\ell$-groups, denoted by $ \ten_{a}$, was defined in \cite{BVanR}. Note then $ \ten_{a}$ is uniquely defined, up to isomorphism, by a universal property with respect to Archimedean structures.

The  tensor product of MV-algebras was defined in \cite{Mun} in both the  non-semisimple  and the semisimple case. 
For two MV-algebras $A$ and $B$, their the tensor product is the MV-algebras $A\otimes_{MV} B$ together with a universal  bimorphism  $\beta_{A,B}:A \times B \rightarrow A \ten_{MV}B$. A {\em bimorphism} is a bilinear function that is $\vee$-preserving and $\wedge$-preserving in each component. The universal property satisfied by $\beta_{A,B}$ is the following: for any MV-algebra $C$ and for any bimorphism $\beta\colon A\times B\rightarrow C$, there is a unique homomorphism of MV-algebras $\omega \colon A\otimes_{MV}B\rightarrow C$ such that $\omega \circ \beta_{A,B}=\beta$. For $a\in A$ and $b\in B$ we denote $a\ten_{MV}b=\beta_{A,B}(a,b)$. As expected,  $A\ten_{MV}B $ is generated by  $\beta_{A,B}(A\times B)$.

Since the class of semisimple MV-algebras is not closed with respect to tensor products,  the tensor product of semisimple MV-algebras is defined in \cite{Mun}  by \[ A\ten B = \quot{(A\ten_{MV}B) }{ Rad(A\ten_{MV}B)},\]
where $Rad(A\ten_{MV}B)$ is the intersection of the maximal ideals of  $A\ten_{MV}B$ and $\ten$ satisfies the same universal property of $\ten_{MV}$ with respect to semisimple MV-algebras.  Recalling that semisimple MV-algebras are isomorphic to subalgebras of continuous functions, we get the following crucial functional representation of $\ten$.
\bthm\cite[Theorem 4.3]{Mun} \label{teo:reprSS}
Let $A,B$ be semisimple MV-algebras and let $X$ and $Y$ be  compact Hausdorff spaces such that $A\subseteq C(X)$ and $B\subseteq C(Y)$. Then $A\ten B$ is an MV-subalgebra of $C(X\times Y)$. Moreover,
\[ A\ten B =< \pi (a, b) \mid a\in A \subseteq C(X),\ b \in B \subseteq C(Y)>_{MV} \subseteq C(X\times Y)\]
where $\pi(a,b)$ is the usual product between functions.
\ethm

\noindent Further properties of the semisimple tensor product of MV-algebras are proved in \cite{LLTP1}. We summarize some of these results in the following theorems.

\bprop
Let $A$ and $B$ be semisimple MV-algebras. Then $A\ten B\simeq B\ten A$, that is, the tensor product of semisimple MV-algebras is commutative.
\eprop
\begin{proof}
It is a straightforward consequence of the functional representation from \cite{Mun}, \Cref{teo:reprSS}.
\end{proof}

\bthm  \label{teo:SEP}
The following hold.
\begin{enumerate}
\item If $A$ is a Riesz MV-algebra and $B$ is a semisimple MV-algebra, $A\otimes B$ is a Riesz MV-algebra.
\item If $A$ and $B$ are unital and semisimple PMV-algebras, $A\otimes B$ is a unital and semisimple PMV-algebra.
\item If $A$ is a unital and semisimple \textit{f}MV-algebra and $P$ is a unital and semisimple PMV-algebra, $A\otimes P$ is a unital and semisimple \textit{f}MV-algebra.
\end{enumerate}
\ethm

\bthm \label{teo:commTP}
If $(G_A, u_A)$, $(G_B, u_B)$ are Archimedean  $\ell u$-groups and  $A$, $B$ are semisimple MV-algebras such that  $A\simeq \Gamma (G_A, u_A)$ and  $B\simeq \Gamma (G_B, u_B)$ then $A\ten B\simeq  \Gamma (G_A \ten_{a}G_B, u_A \ten _{a}u_B)$.
\ethm

Finally, in \cite{LLTP1} the following functors are defined.
\begin{itemize}
\item ${\mathcal T}_\ten:\mathbf{MV_{ss}}\to \mathbf{RMV_{ss}}$ is defined by 
\begin{itemize}
\item ${\mathcal T}_\ten(B)=[0,1]\ten B$, which is a semisimple Riesz MV-algebra by Theorem \ref{teo:SEP};

\item for any homomorphism of MV-algebras $f:A\to B$, ${\mathcal T}_\ten(f) =\widetilde{f}$ where $\widetilde{f}:[0,1]\ten A\to [0,1]\ten B$ is the unique Riesz MV-algebra homomorphism such that $\widetilde{f}\circ \iota_A=\iota_B\circ f$, which exists by \cite[Corollary 3.1]{LLTP1}.
\end{itemize}
\item $\mathcal{F}_{\ten}: \mathbf{uPMV_{ss}} \rightarrow \mathbf{ufMV_{ss}}$ is defined by
\begin{itemize}
\item $\mathcal{F}_{\ten}(P)=[0,1]\ten P$, which is  a unital, commutative and semisimple \textit{f}MV-algebra by Theorem \ref{teo:SEP};
\item for any homomorphism of PMV-algebras $h:P_1\rightarrow P_2$, $\mathcal{F}_\ten(h)=h^{\sharp}$  is the homomorphism of  \textit{f}MV-algebras defined in \cite[Proposition 4.3]{LLTP1}, that satisfies the condition $h^{\sharp} \circ \iota_1 = \iota_2 \circ h$.
\end{itemize}
\item From $\mathbf{ufMV_{ss}}$ to $\mathbf{uPMV_{ss}}$ and from $\mathbf{RMV_{ss}}$ to $\mathbf{MV_{ss}}$ we have the usual forgetful functor $\mathcal{U}_{\mathbb{R}}$.\\
\end{itemize}

\bthm \cite{LLTP1}
($\mathcal{T}_\ten, \mathcal{U}_{\mathbb R})$  and ($\mathcal{F}_{\ten}, \mathcal{U}_{\mathbb R}$) are two pairs of adjoint functors.
\ethm

By categorical equivalence, the adjunctions $(\mathcal{T}_\ten , \mathcal{U}_{\mathbb{R}})$ and $(\mathcal{F}_\ten , \mathcal{U}_{\mathbb{R}})$ are naturally transfered to lattice-ordered groups and rings, Riesz spaces and algebras: we obtain an adjunction between $\mathbf{auG_a}$ and $\mathbf{uRS_a}$, that we shall denote by  $(\mathcal{T}_{\ten a} , \mathcal{U}_{\ell \mathbb{R}})$ and an adjunction between $\mathbf{uR_a}$ and $\mathbf{fuAlg_a}$, denoted by $(\mathcal{F}_{\ten a} , \mathcal{U}_{\ell \mathbb{R}})$.

Figure 1 merges together the results present in literature. We remark that the adjunctions represented by the horizontal lines were proved by the authors in \cite{LLTP1} and represent the stepping stone of the present investigation.
\begin{center}
\begin{tikzpicture}
\node (A) {$\mathbf{auG_{a}}$};
\node (B) [below=1.5cm of A, right of=A]{$\mathbf{MV_{ss}}$};
\node (C) [right of=B] {$\mathbf{RMV_{ss}}$};
\node (D) [above=1.5cm of C, right of=C] {$\mathbf{uRS_{a}}$};
\node (E) [below of=B] {$\mathbf{uPMV_{ss}}$};
\node (F) [right of=E] {$\mathbf{ufMV_{ss}}$};
\node (G) [below=1.5cm of E, left of=E] {$\mathbf{uR_{a}}$};
\node (H) [below=1.5cm of F, right of=F] {$\mathbf{ufAlg_{a}}$};
\draw[<->] (A) to node {$\Gamma$} (B);
\draw[->] (D) to node {$\mathcal{U}_{(\ell\mathbb R)}$} (A);
\draw[->] (C) to node {$\mathcal{U}_{\mathbb{R}}$} (B);
\draw[<->] (C) to node {$\Gamma_{\mathbb R}$} (D);
\draw[->] (E) to node {$\mathcal{U}_{(\cdot)}$} (B);
\draw[->] (F) to node [swap]{$\mathcal{U}_{\mathbb R}$} (E);
\draw[->] (F) to node [swap]{$\mathcal{U}_{(\cdot)}$} (C);
\draw[->] (G) to node {$\mathcal{U}_{(\cdot \ell)}$} (A);
\draw[->] (H) to node [swap]{$\mathcal{U}_{(\cdot \ell)}$} (D);
\draw[->] (H) to node [swap]{$\mathcal{U}_{(\ell \mathbb R)}$} (G);
\draw[<->] (G) to node {$\Gamma_{(\cdot)}$} (E);
\draw[<->] (F) to node {$\Gamma_{f}$} (H);

\draw[->, bend right] (E) to node [swap]{{$\mathcal{F}_{\ten}$}} (F);
\draw[->, bend left] (B) to node {{${\mathcal T}_\ten$}} (C);

\draw[->, bend right=12] (G) to node [swap]{{$\mathcal{F}_{\ten a}$}} (H);
\draw[->, bend left=12] (A) to node {{${\mathcal T}_{\ten a}$}} (D);

\node (S)  at ($(G)!0.5!(H)$) { };
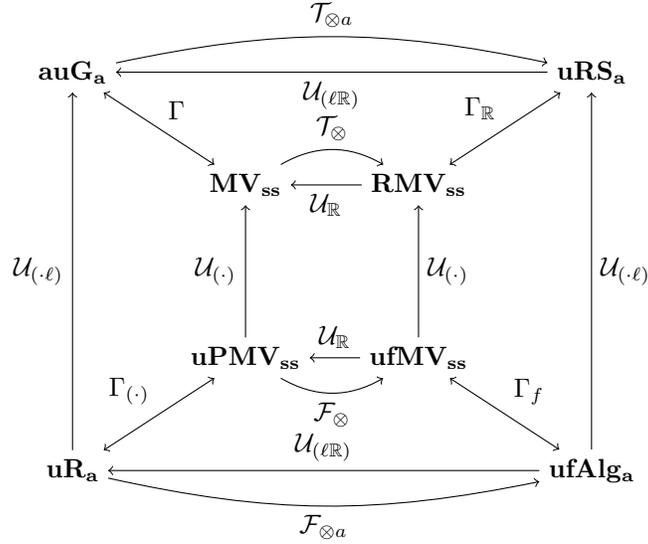
\node [below=1cm, align=flush center] at (S){ Figure 1. Relations between categories.};
\end{tikzpicture}
\end{center}

\subsection{The non-semisimple case}\label{sec:SEP}
As discussed in the previous section, the tensor product of MV-algebras has been defined in the non-semisimple case as well, although it does not enjoy the very convenient functional representation presented in \Cref{teo:reprSS}. In the same way, the tensor product of $\ell u$-groups is given in the non-semisimple case, and a variant of \Cref{teo:commTP} is proved in \cite{LLTP1}.

The approach we followed in \cite{LLTP1}, and that we are continuing pursuing here, makes fundamental use of the \emph{scalar extension property}. Such a property is easily obtained in the case of non-ordered structures, but we could not prove it for non-semisimple MV-algebras, the problems being the presence of infinitesimal elements, the lattice-order and the fact that, in order to define the scalar multiplication on the tensor product, one has to do it by universal property. Altogether,  these fact entails that, to obtain the desired scalar extension property, one needs to prove that the sum of two peculiar MV-homomorphisms is an MV-homomorphism,  but this is not always true. We also remark that the same property of sums of homomorphisms will be needed subsequently to define the tensor PMV-algebra of an MV-algebra. Moreover, the same problems appear in the setting of $\ell$-groups: in \cite{Mart}, the proof is a sketch that we were not be able to complete. It is worth mentioning that the statement of the scalar extension property for $\ell$-groups can be found in other references, but again (and to the best of our knowledge) is given without proof (or with a sketch of it) in each instance. Whence, in this paper we restrict our attention to semisimple structures. In other fields this may seem a set-back, but it is more than adequate for our logic-oriented approach. Indeed, looking back to logic one can see that the most significant algebras to work with are free algebras and finitely presented algebras.  Such algebras correspond to theories (in the case of finitely presented algebras, to finitely axiomatizable theories) in the logic associated to the class of algebras, and they are both semisimple. Moreover, the \emph{standard model} for all logics considered is the simple algebra $[0,1]$ endowed with appropriate operations.

\section{The semisimple tensor PMV-algebra of a semisimple MV-algebra} \label{sec:TPMV}
In this section we provide the algebraic tools needed to complete Figure 1 with adjoints for all forgetful functors.  The key ingredient will be an ``MV-algebraic'' version of the classical construction of the \textit{tensor algebra}. 

We start by proving that the semisimple tensor product of MV-algebras is associative. We recall again that any semisimple MV-algebra is isomorphic to a subalgebra of $C(X)$ \--- the MV-algebra of continuous and $[0,1]$-valued functions defined on $X$ \--- for a suitable compact Hausdorff space $X$, namely, the space of its maximal ideals.

\bprop \label{pro:ass01}
Let $A,B,C$ be semisimple MV-algebras and let $X,Y,Z$ be the spaces such that $A \subseteq C(X)$, $B \subseteq C(Y)$, $C \subseteq C(Z)$. Then $ A\ten (B\ten C) =   (A\ten B)\ten C= \left\langle a \cdot b \cdot c\mid a\in A, b\in B, c\in C \right\rangle \subseteq C(X\times Y \times Z).$
\eprop
\begin{proof}
Let $M$ be the MV-subalgebra of $C(X\times Y\times Z)$ generated by $a\cdot b\cdot c$, where $\cdot $ is the usual product between functions. By Theorem \ref{teo:reprSS}, $(A\ten B)\ten C =\left\langle f \cdot c \mid f\in A\ten B, c\in C \right\rangle $, the MV-algebra generated by the product of $f\in A\ten B$ and $c\in C$. We want to prove that $ \left\langle f \cdot c \mid f\in A\ten B, c\in C \right\rangle = M.$
Trivially $M\subseteq \left\langle f \cdot c \mid f\in A\ten B, c\in C \right\rangle$. We prove the other inclusion by induction on the construction of $f\in A\ten B$.

(i) If $f=a \ten b = a \cdot b$, then it is trivial: $(a \cdot b )\cdot c =a\cdot b\cdot c \in M$.

(ii) Let $f$ be in $A\ten B$ such that $f\cdot c\in M$. Then $f^* \cdot c=(\mathbf{1}-f)\cdot c= c - f\cdot c \in M $
 by induction hypothesis and the fact that $\mathbf{0}\le c-f\cdot c\le \mathbf{1}$ in $C(X\times Y\times Z)$.

(iii) Let $f=f_1 \oplus f_2$ be in $A\ten B$ such that $f_1 \cdot c$ and $f_2 \cdot c$ belongs to $M$. Since we deal with subalgebras of continuous functions, 
it is easily seen that the product distributes over $\wedge$, and therefore we have

 $f \cdot c = (f_1 \oplus f_2) \cdot c = (f_1 + (f_1^* \wedge f_2)) \cdot c = f_1 \cdot c + (f_1^* \cdot c \wedge f_2 \cdot c) $\\ which belongs to $M$ by induction hypothesis and (ii).\\
Therefore $M= (A\ten B) \ten C$. In the same way we prove that $M=A\ten (B\ten C)$ and the claim is settled.
\end{proof}

\bcor \label{cor:ass01}
The tensor product of semisimple MV-algebras is associative.
\ecor

\bfact The same result for $\ten_{MV}$ would need an analogous of Theorem \ref{teo:reprSS}, which relies on the functional representation of semisimple MV-algebras and H\"older's theorem for Archimedean unital $\ell$-groups. A possible solution for the non-semisimple case is an open problem.
\efact
 
\bdfn 
Let $A$ be a semisimple MV-algebra and let $X$ be the compact Hausdorff space such that $A\subseteq C(X)$. We define:
\[T^1(A)=A, \quad T^n(A)=T^{n-1}(A)\ten A,\]
where $\ten$ is the semisimple tensor product. By Proposition \ref{pro:ass01}, 
\[T^n(A) = \left\langle f_1 \cdot \ldots \cdot f_n \mid f_i \in A, i=1 \ldots n \right\rangle \subseteq C(X^n),\]
and  $\overbrace{1\ten \ldots \ten 1}^{n}= \overbrace{1\cdot \ldots \cdot 1}^{n}$ is the top element of $T^n(A)$ for every $n$.\\
For any $n, m \in \mathbb{N}$  with $n\le m$  we define
\begin{itemize}
\item[] $\epsilon_{n,n}$ as the identity homomorphism on $T^n(A)$;
\item[] $ \epsilon_{n,m}: T^n(A) \rightarrow T^m(A)$, by $\epsilon_{n,m}(\mathbf{x})= \mathbf{x}\ten (1\ten \ldots \ten 1)$,
\end{itemize}
where by associativity $T^m(A) \simeq T^n(A) \ten T^{m-n}(A)$. By \cite[Proposition 2.1]{LLTP1}, $\epsilon_{n,m}$ is the embedding in the semisimple tensor product and $\epsilon_{m,k}\circ \epsilon_{n,m}=\epsilon_{n,k}$. We remark that any $T^n(A)$ is semisimple by construction.\\
$(T^n(A), \epsilon_{n,m})$ is a direct system, therefore we consider the disjoint union 
\[\bigsqcup_{n\in \mathbb{N}}T^n(A),\]
and we define an equivalence relation on it by
\[(x,n)\sim (y,m)\ \Leftrightarrow \text{ there exists }k \ge n,m \text{ such that }\epsilon_{n,k}(x)=\epsilon_{m,k}(y).\]
The quotient MV-algebra $T(A) = \quot{ \bigsqcup_{n \in \mathbb{N}}T^n(A)}{\sim }$
is the direct limit of the direct system, and $\epsilon _{n, A} :T^n(A) \rightarrow T(A)$ is the canonical morphism that maps each element in its equivalence class. When there is no danger of confusion, we will denote $\epsilon_{n,A}$ simply by $\epsilon_n$. 

We call $T(A)$ the \textit{Tensor} PMV-\textit{algebra} of the MV-algebra $A$.
\edfn

\blem
The algebra $T(A)$ is a semisimple MV-algebra.
\elem
\begin{proof}
Suppose that there exist an infinitesimal element $\mathbf{x} \in T(A)$. It follows that $n \mathbf{x}\le \mathbf{x}^*$ for any $n\in \mathbb{N}$, therefore $n \mathbf{x} \odot \mathbf{x}=\mathbf{0}$ for any $n\in \mathbb{N}$. This comes to the existence of naturals $m,l,k$ such that $\mathbf{x}$ is the equivalence class of $(x, m)$, $\mathbf{0}$ is the equivalence class of $(0, l)$ and $n\epsilon_{m,k}(x)\odot \epsilon_{m,k}(x)=\epsilon_{l,k}(0)$. This entails that $\epsilon_{m,k}(x)$ is infinitesimal in the semisimple MV-algebra $T^k(A)$, a contradiction.
\end{proof}

\bfact The above construction is the definition of the direct limit of a direct system in category theory, specialized to our framework. Therefore, it is well known that  $\epsilon_m \circ \epsilon_{n,m}= \epsilon_n$ for any $n\le m$ and that the limit $(T(A),\epsilon _n)$ enjoys the universal property that makes the diagram in Figure 2 commutative. 
\begin{center}
\begin{tikzpicture}
  \node (X) {$T^n(A)$};
  \node (T) [right of=X] {$T^m(A)$};
  \node (Y) [right of=T] {$T^k(A)$ };
  \node (Z) [below=1cm of T] {$T(A)$};
  \draw[->] (X) to node {$\epsilon _{n,m}$} (T);
  \draw[->] (X) to node [swap]{$\epsilon _n$} (Z);
  \draw[->] (T) to node [swap]{$\epsilon _m$} (Z);
  \draw[->] (T) to node {$\epsilon_{m,k}$} (Y);
  \draw[->] (Y) to node {$\epsilon_k$} (Z);
  \node [below=0.5cm, align=flush center] at (Z){ Figure 2. Direct limit.};
 \end{tikzpicture}
\end{center}
\efact

\bnot
For any $\mathbf{a}\in T^n(A)$ and any $ \mathbf{b}\in T^m(A) $ in order to avoid confusion, we denote the bimorphism $\pi$ from  Theorem \ref{teo:reprSS} by
\begin{itemize}
\item[] $\gamma_{n,m}: T^n(A) \times T^m(A) \rightarrow T^{n+m}(A) \subseteq C(X^{n+m})$,
\item[] $\gamma_{n,m}(\mathbf{a}, \mathbf{b})(x_1, \ldots, x_n, y_1, \ldots, y_m)= \mathbf{a}(x_1, \ldots , x_n)\mathbf{b}(y_1,\ldots , y_m)$.
\end{itemize}
\enot

The following lemma collects some technical properties of the maps $\epsilon_{n,m}$ and $\gamma_{n,m}$. All proofs rely on the functional representation of the algebras $T^k(A)$, for any semisimple MV-algebra $A$.

\blem \label{lem:propTPMV}
For any $n, m , k \in \mathbb{N}$, the following hold:
\begin{enumerate}[label=(\roman*)]
\item  $\gamma_{n,m}(\mathbf{a}, \mathbf{1}_m)= \epsilon_{n,n+m}(\mathbf{a})$, with $\mathbf{a}\in T^n(A)$ and $\mathbf{1}_m$ top element in $T^m(A)$, that is unit function in $C(X^m)$.
\item $\epsilon_{n+m}=\epsilon_{m+n}$ and $\epsilon_{n+(m+l)}=\epsilon_{(n+m)+l}$;
\item $\gamma_{n,m}(\mathbf{a}, \mathbf{b})=\gamma_{m,n}(\mathbf{b}, \mathbf{a})$, for any $\mathbf{a}\in T^n(A)$ and $\mathbf{b}\in T^m(A)$;
\item If $n\le m$,  $\gamma_{n, m+k}(\mathbf{a}, \gamma_{m,k}(\mathbf{b}, \mathbf{c}))=\gamma_{n+m, k}(\gamma_{n,m}(\mathbf{a}, \mathbf{b}), \mathbf{c})$, for any $\mathbf{a}\in T^n(A)$, $\mathbf{b}\in T^m(A)$ and $\mathbf{c}\in T^k(A)$;
\item If $n\le m$, $\gamma_{m,k}(\epsilon_{n,m}(\mathbf{a}), \mathbf{b})=\epsilon_{n+k, m+k}(\gamma_{n,k}(\mathbf{a}, \mathbf{b}))$.
\end{enumerate}
\elem
\begin{proof}

(i) It is straightforward by definition of all $\gamma$-maps and $\epsilon$-maps.

(ii) It is straightforward by Proposition \ref{pro:ass01} and the universal property of the direct limit.

(iii) We recall that any $T^l(A)$ is a subalgebra of $C(X^l)$. Moreover, $\gamma_{n,m}(\mathbf{a}, \mathbf{b})= \mathbf{a}\cdot \mathbf{b} \in C(X^{n+m})$ and $\gamma_{m,n}(\mathbf{b}, \mathbf{a})=\mathbf{b}\cdot \mathbf{a}\in C(X^{m+n})$. Since $X^{n+m}\simeq X^{m+n}$, the conclusion follows by the commutativity of the product of functions.

(iv) By definition, $\gamma_{n, m+k}(\mathbf{a}, \gamma_{m,k}(\mathbf{b}, \mathbf{c}))= \mathbf{a}\cdot (\mathbf{b}\cdot \mathbf{c}) \in C(X^{n+(m+k)})$ and $\gamma_{n+m, k}(\gamma_{n,m}(\mathbf{a}, \mathbf{b}), \mathbf{c})= (\mathbf{a}\cdot \mathbf{b})\cdot \mathbf{c} \in C(X^{(n+m)+k)})$. Since $X^{n+(m+k)}\simeq X^{(n+m)+k}$, the conclusion follows by the associativity of the product of functions.

(v) By (1), (3) and (4) we have
\begin{eqnarray}
&&\gamma_{m,k}(\epsilon_{n,m}(\mathbf{a}), \mathbf{b})= \gamma_{(m-n)+n,k}(\gamma_{n,m-n}(\mathbf{a}, \mathbf{1}_{m-n}), \mathbf{b})=\nonumber \\
&=&\gamma_{(m-n)+n,k}(\gamma_{m-n,n}(\mathbf{1}_{m-n}, \mathbf{a}), \mathbf{b})=\gamma_{m-n, n+k}(\mathbf{1}_{m-n}, \gamma_{n,k}(\mathbf{a}, \mathbf{b}))=\nonumber \\
&=&\gamma_{ n+k, m-n}(\gamma_{n,k}(\mathbf{a}, \mathbf{b}), \mathbf{1}_{m-n} )= \epsilon_{n+k, m+k}(\gamma_{n,k}(\mathbf{a}, \mathbf{b})).\nonumber
\end{eqnarray}
\end{proof}
\bprop \label{pro:tenPMV}
For any semisimple MV-algebra $A$, $T(A)$ is a semisimple and unital PMV-algebra.
\eprop 
\begin{proof}
We define the product as follows. For any $\mathbf{x}, \mathbf{y}\in T(A)$ there exist $n, m \in \mathbf{N}$ such that $\mathbf{x}= \epsilon_n (\mathbf{a})$, with $\mathbf{a}\in T^n(A)$ and $\mathbf{y}= \epsilon_m (\mathbf{b})$, with $\mathbf{b}\in T^m(A)$. Then 

\[ \mathbf{x} \cdot \mathbf{y} = (\epsilon_{n+m} \circ \gamma_{n,m})(\mathbf{a}, \mathbf{b}).\] 
We first need to prove that the operation is well defined. Let $\mathbf{c}\in T^l(A)$ and $\mathbf{d}\in T^k(A)$ be elements such that $(\mathbf{a}, n)\sim (\mathbf{c}, l)$ and $(\mathbf{b}, m)\sim (\mathbf{d}, k)$. This means that we can assume, without loss of generality, $\mathbf{a}=\epsilon_{l,n}(\mathbf{c})$ and $\mathbf{b}=\epsilon_{k,m}(\mathbf{d})$.  Then, applying Lemma \ref{lem:propTPMV}, we get
\begin{eqnarray}
&&\epsilon_{n+m}(\gamma_{n,m}(\mathbf{a},\mathbf{b}))=\epsilon_{n+m}(\gamma_{n,m}(\epsilon_{l,n}(\mathbf{c}), \epsilon_{k,m}(\mathbf{d}))) \nonumber\\
&=& \epsilon_{n+m}(\epsilon_{m+l, m+n}(\gamma_{l,m}(\mathbf{c}, \epsilon_{k,m}(\mathbf{d}))))= \epsilon_{m+l}(\gamma_{m,l}(\epsilon_{k,m}(\mathbf{d}), \mathbf{c})) \nonumber \\
&=& \epsilon_{m+l} (\epsilon_{k+l, m+l}(\gamma_{l,k}(\mathbf{c},\mathbf{d})))= \epsilon_{l+k}(\gamma_{l,k}(\mathbf{c},\mathbf{d})).\nonumber
\end{eqnarray}
To prove that $T(A)$ is a PMV-algebra, let us prove that  the function $(\mathbf{x}, \mathbf{y})\mapsto \mathbf{x}\cdot \mathbf{y}$ is bilinear, that the product is associative and that $\mathbf{1}$ is the unit.\\
To prove bilinearity\footnote{We recall that in order to prove linearity of a map $\beta:A \to B$, one needs to prove that $\beta (a_1\oplus a_2)=\beta (a_1)\oplus \beta (a_2)$ for any $a_1, a_2$ such that $a_1\le a_2^*$. To give a more compact version of the statement, it is possible to define a partial sum $+$ by ``$a_1 + a_2$ is defined iff $a_1\le a_2^*$ and in this case $a_1 +a_2=a_1\oplus a_2$''. This is equivalent to state that $a_1 + a_2$ is defined iff $a_1 \odot a_2=0$, see \cite{CDM}. Thus, linearity of $\beta$ equals to $\beta(a_1 + a_2)=\beta(a_1)+\beta(a_2)$.}, let $\mathbf{x_1}, \mathbf{x_2}, \mathbf{y}$ be elements in $T(A)$ such that 
\begin{eqnarray}
\mathbf{x_1}=\epsilon_n(\mathbf{a_1}) \text{ with } \mathbf{a_1}\in T^n(A),&& \mathbf{x_2}=\epsilon_m(\mathbf{a_2}) \text{ with } \mathbf{a_2}\in T^m(A) \nonumber \\
\mathbf{x_1}+ \mathbf{x_2} \text{ defined, } && \mathbf{y}=\epsilon_k(\mathbf{c}) \text{ with } \mathbf{c}\in T^k(A).\nonumber
\end{eqnarray}
\noindent Without loss of generality, we assume that $n\le m$. Therefore, since each of the $\epsilon$-map is a homomorphism of MV-algebras,
\[\epsilon_n(\mathbf{a_1})+ \epsilon_m(\mathbf{a_2})=\epsilon_m(\epsilon_{n,m}(\mathbf{a_1}))+\epsilon_m(\mathbf{a_2})=\epsilon_m(\epsilon_{n,m}(\mathbf{a_1})+\mathbf{a_2}).\]
Thus, by definition of the product on $T(A)$, 

\[(\mathbf{x_1}+\mathbf{x_2})\cdot \mathbf{y}= \epsilon_{m+k}(\gamma_{m,k}(\epsilon_{n,m}(\mathbf{a_1})+\mathbf{a_2}, \mathbf{b}))\]
and since $\gamma_{m,k}$ is a bimorphism,
\[\gamma_{m,k}(\epsilon_{n,m}(\mathbf{a_1})+\mathbf{a_2}, \mathbf{b})=\gamma_{m,k}(\epsilon_{n,m}(\mathbf{a_1}), \mathbf{b})+\gamma_{m,k}(\mathbf{a_2}, \mathbf{b}).\]
By Lemma \ref{lem:propTPMV} (5) we have $\gamma_{m,k}(\epsilon_{n,m}(\mathbf{a_1}), \mathbf{b})=\epsilon_{n+k, m+k}(\gamma_{n,k}(\mathbf{a_1}, \mathbf{b}))$ and
\begin{align*}
(\mathbf{x_1}+\mathbf{x_2})\cdot \mathbf{y}=& \epsilon_{m+k}(\epsilon_{n+k, m+k}(\gamma_{n,k}(\mathbf{a_1}, \mathbf{b}))) + \epsilon_{m+k}(\gamma_{m,k}(\mathbf{a_2}, \mathbf{b}))=\\
=& \epsilon_{n+k}(\gamma_{n,k}(\mathbf{a_1}, \mathbf{b}))+ \mathbf{x_2}\cdot \mathbf{y}= \mathbf{x_1}\cdot \mathbf{y}+ \mathbf{x_2}\cdot \mathbf{y}.
\end{align*}
\noindent One can prove in the same way that  $\mathbf{y}\cdot (\mathbf{x_1}+\mathbf{x_2})= \mathbf{y}\cdot \mathbf{x_1}+\mathbf{y}\cdot \mathbf{x_2}$.

Associativity follows directly from Lemma \ref{lem:propTPMV} (2) and (4).

Finally, for any $k\in \mathbf{N}$ we denote by $\mathbf{1}$ and $\mathbf{1}_k$ the top elements of $T(A)$ and $T^k(A)$ respectively.  It follows $\epsilon_k (\mathbf{1}_k)=\mathbf{1}$ for any $k\in \mathbb{N}$.

Let $\mathbf{x}\in T(A)$, such that $\mathbf{x}=\epsilon_n(\mathbf{a})$ with $\mathbf{a}\in T^n(A)$ and let $m$ be a positive integer such that $\mathbf{1}=\epsilon_m(\mathbf{1}_m)$. We have 
\[\mathbf{x}\cdot \mathbf{1}=\epsilon_{n+m}(\gamma_{n,m}(\mathbf{a}, \mathbf{1}_m))=\epsilon_{n+m}(\epsilon _{n, n+m}(\mathbf{a}))=\epsilon_{n} (\mathbf{a}) =\mathbf{x}.\]
The proof of the equality $\mathbf{1}\cdot \mathbf{x}=\mathbf{x}$ follows from Lemma \ref{lem:propTPMV} (1) and (3), then $T(A)$ is a unital PMV-algebra whose MV-algebra reduct is semisimple. This entails that $T(A)$ is a unital and semisimple PMV-algebra.
\end{proof}

\blem \label{lem:emb1}
Any map $\epsilon_n : T^n(A) \rightarrow T(A)$ is an embedding.
\elem
\begin{proof}
It is straightforward by \cite[\S\ 21, Lemma 2]{Gr}.
\end{proof}

\bthm \label{teo:tenUnProp}
Let $A$ be a semisimple MV-algebra. For any semisimple and unital PMV-algebra $P$ and for any homomorphism of MV-algebras $f: A\rightarrow \mathcal{U}_{(\cdot)}(P)$ there exists a homomorphism of PMV-algebras $f^{\sharp}:T(A) \rightarrow P$ such that $\widetilde{f} \circ \epsilon_{1,A} = f$.
\ethm
\begin{proof}
Using the universal property of the tensor product, let us define the following family of homomorphisms $\{\widetilde{\lambda_n}\}_{n\in \mathbb{N}}$.
\begin{enumerate}[label=(\roman*)]
\item For $n=1$, $\lambda_1= \widetilde{\lambda}_1=f$.
\item For $n=2$, we define $\lambda_2: A\times A \rightarrow P$ to be the function such that $\lambda_2 (a_1,a_2)=f(a_1)\cdot f(a_2)$.  Since $P$ is a unital PMV-algebra, $\lambda_2$ is a bimorphism and $\lambda_2 (1_A,1_A)=f(1_A)\cdot f(1_A)=1_P\cdot 1_P=1_P$. Whence, there exists a homomorphism of MV-algebras $\widetilde{\lambda_2}: A\ten A \rightarrow A$ such that $\widetilde{\lambda}_2(a_1\ten a_2)=f(a_1)\cdot f(a_2)$.
\item for any $n\in \mathbb{N}$, $\lambda_n: T^{n-1}(A)\times A \rightarrow P$, $\lambda_n (\mathbf{x},a_n)=\widetilde{\lambda}_{n-1}(\mathbf{x})\cdot f(a_n)$. Thus, $\widetilde{\lambda}_n: T^n(A)\rightarrow P$ is the homomorphism such that $\widetilde{\lambda}_n(a_1\ten \ldots \ten a_n)=f(a_1)\cdot \ldots \cdot f(a_n).$
\end{enumerate}
It is easily seen that $\widetilde{\lambda}_m \circ \epsilon_{n,m}= \widetilde{\lambda}_n$ for any $n\le m$, since the two homomorphisms coincide on generators.

Thus, we have the following situation
\begin{center}
\begin{tikzpicture}
  \node (X) {$T^n(A)$};
  \node (T) [right of=X] {};
  \node (Y) [right of=T] {$T^m(A)$ };
  \node (Z) [below=1cm of T] {$T(A)$};
  \node (W) [below=1 of Z] {$P$};
  \draw[->] (X) to node {$\epsilon _{n,m}$} (Y);
  \draw[->] (X) to node {$\epsilon _n$} (Z);
  \draw[->] (Y) to node [swap]{$\epsilon _m$} (Z);
  \draw[->] (X) to node [swap]{$\widetilde{\lambda}_n$} (W);
  \draw[->] (Y) to node {$\widetilde{\lambda}_m$} (W);
  \draw[->][dashed] (Z) to node [swap]{ } (W);
   \node [below=0.5cm, align=flush center] at (W){ Figure 3. Universal property.};
 \end{tikzpicture}
\end{center}
Since $(T(A), \epsilon _n ) $ is the direct limit for $( T^n(A), \epsilon _{n,m} ) $, there exists a homomorphism of MV-algebras $\widetilde{f}: T(A) \rightarrow P$ such that $ \widetilde{f}\circ \epsilon _n = \widetilde{\lambda}_n$, and then $\widetilde{f}\circ \epsilon_1 = \widetilde{\lambda_1}=f$.\\
Finally, the fact that $\widetilde{f}$ is an homomorphism of PMV-algebras is a direct consequence of \cite[Theorem 1.1]{BH}, taking into account the fact that unital and semisimple PMV-algebras correspond unital and Archimedean $f$-rings, and that such rings are semiprime (that is, without non-trivial nilpotent elements). 
\end{proof}
\bfact \label{rem:univProp}
A direct proof of the fact that $\widetilde{f}$ is an homomorphism of PMV-algebras can be found in \cite{LPhD}.
\efact

\bcor \label{pro:tenhom}
Let $A$, $B$ be semisimple MV-algebras and $h:A\rightarrow B$ be a homomorphism of MV-algebras. Then there exists a unique homomorphism of PMV-algebras $h^{\sharp}:T(A)\rightarrow T(B)$ such that $h^{\sharp}\circ \epsilon_{1,A} = \epsilon_{1,B} \circ h$.
\ecor
\begin{proof}
We apply \Cref{teo:tenUnProp} with $f=\epsilon_{1,B} \circ h$.
\end{proof}

We now apply the construction of the Tensor PMV-algebra to a Riesz MV-algebra.
\bthm \label{thm:tenfMV}
Let $A$ be a semisimple Riesz MV-algebra. Then $T(A)$ is a unital and semisimple \textit{f}MV-algebra.
\ethm
\begin{proof}
We recall that by Theorem \ref{teo:SEP}, any $T^n(A)$ is a Riesz MV-algebra.\\
Let $\mathbf{x}\in T(A)$, then there exist $n\in \mathbb{N}$ and $\mathbf{a}\in T^n(A)$ such that $\mathbf{x}=\epsilon_n(\mathbf{a})$. 
We define the scalar product on the direct limit as 
\begin{center}
$ \alpha \mathbf{x}= \epsilon_n(\alpha \mathbf{a}),$ for any $\alpha \in [0,1]$.
\end{center}
It is easily seen that the operation is well defined, using the fact that each $\epsilon_{n,m}$ is an homomorphism of Riesz MV-algebras\footnote{This follows from the fact that any MV-homomorphism between Riesz MV-algebras preserve the scalar product, see \cite[Theorem 2]{LeuRMV}. We remark that the aforementioned corollary actually requires for the codomain of the homomorphism to be a semisimple algebra \cite[Proposition 11.53]{ce}, which the case of interest for us.}. 

Take $\mathbf{x}=\epsilon_n(\mathbf{a})$, $\mathbf{y}=\epsilon_m(\mathbf{b})$ and, without loss of generality, assume $n\le m$. If the partial sum is defined we get
\[\mathbf{x}+\mathbf{y}= \epsilon_n(\mathbf{a})+\epsilon_m(\mathbf{b}) = \epsilon_m (\epsilon_{n,m}(\mathbf{a}))+\epsilon_m(\mathbf{b})=\epsilon_m (\epsilon_{n,m}(\mathbf{a})+\mathbf{b}).\]
As any $T^n(A)$ a Riesz MV-algebra, and since any $\epsilon_{n,m}$ a homomorphism of Riesz MV-algebras, we infer:
\begin{enumerate}[label=(\roman*)]
\item $\alpha (\mathbf{x}+\mathbf{y})=\epsilon_m (\alpha (\epsilon_{n,m}(\mathbf{a})+ \mathbf{b}))= \epsilon_m ( \epsilon_{n,m}(\alpha \mathbf{a}))+ \epsilon_m(\alpha \mathbf{b})= \epsilon_n (\alpha \mathbf{a})+ \epsilon_m(\alpha \mathbf{b})=  \alpha \mathbf{x}+ \alpha \mathbf{y}$.
\item $(\alpha + \beta)\mathbf{x}= \epsilon_n((\alpha + \beta)\mathbf{a})=\epsilon_n (\alpha \mathbf{a}+\beta \mathbf{a})= \epsilon_n(\alpha \mathbf{a})+\epsilon_n (\beta \mathbf{a})= \alpha \mathbf{x}+\beta \mathbf{x}.$
\item $(\alpha \cdot \beta)\mathbf{x}= \epsilon_n((\alpha \cdot \beta)\mathbf{a})=\epsilon_n (\alpha (\beta \mathbf{a}))= \alpha \mathbf{y}$, with $\mathbf{y}=\epsilon_n (\beta \mathbf{a})$ and $\mathbf{y}=\epsilon_n(\beta \mathbf{a})= \beta \mathbf{x}.$
\item $1\mathbf{x}= \epsilon_n (1\mathbf{a})= \epsilon_n (\mathbf{a})=\mathbf{x}$.
\end{enumerate}

Hence, $T(A)$ is a unital PMV-algebra and a Riesz MV-algebra. Using the fact that each $T^n(A)$ is an algebra of fuctions, we infer the associativity-like required in the definition of a \textit{f}MV-algebras:
\begin{eqnarray}
&&\alpha (\mathbf{x}\cdot \mathbf{y})= \epsilon_{n+m}(\alpha \gamma_{n,n}(\mathbf{a},\mathbf{b}))=\epsilon_{n+m}(\alpha (\mathbf{a}\cdot \mathbf{b}))=\epsilon_{n+m}((\alpha \mathbf{a})\cdot \mathbf{b}))=\nonumber \\
&=&\epsilon_{n+m}( \gamma_{n,m}(\alpha \mathbf{a},\mathbf{b})) =(\alpha \mathbf{x})\cdot \mathbf{y}.\nonumber
\end{eqnarray}
In the same way we prove that $\alpha (\mathbf{x}\cdot \mathbf{y})=\mathbf{x}\cdot (\alpha \mathbf{y})$ and $T(A)$ is a unital and semisimple \textit{f}MV-algebra.
\end{proof}

\bthm
If $A$ is a unital and semisimple PMV-algebra, then $A\simeq T(A)$.
\ethm
\begin{proof}
We first remark that, being $A$ a PMV-algebra, for any $n\in \mathbb{N}$,
\[T^n(A)=\left\langle  f_1 \cdot \ldots \cdot f_n \mid f_i \in A \subseteq C(X)\right\rangle  \subseteq A=T^1(A).\]
Each $T^n(A)$ is therefore an MV-subalgebra of $A$ and in general it is not a PMV-subalgebra. Moreover, since $a_1\cdot \ldots \cdot a_n=a_1\cdot \ldots \cdot a_n\cdot\stackrel{n-1}{\overbrace{1_A\cdot \ldots \cdot 1_A}}\in T^n(A)$, $\epsilon_{n,n}=\epsilon_{1,n}$ for any $n\in \mathbb{N}$.

Let us prove that $\epsilon_1:A\to T(A)$ gives the desired isomorphism.  Indeed, let $\mathbf{y}$ be an element of $T(A)$, with $\mathbf{y}= \epsilon_n (\mathbf{a})$, for some $n\in \mathbb{N}$ and $\mathbf{a}\in T^n(A) \subseteq A$. Then $\mathbf{a}=\epsilon_{1,n}(\mathbf{a})$ and $\mathbf{y}= \epsilon_n (\mathbf{a})= \epsilon_n (\epsilon_{1,n}(\mathbf{a}))=\epsilon_1(\mathbf{a})$, that is, $\epsilon_1$ is surjective.  We prove directly that $\epsilon_1$ is a homomorphism of PMV-algebras. By definition, $\epsilon_1(\mathbf{a}_1 \cdot \mathbf{a}_2)=\epsilon_2(\mathbf{a}_1 \cdot \mathbf{a}_2)=\epsilon_2(\gamma_{1,1}(\mathbf{a}_1, \mathbf{a}_2))=\epsilon_1(\mathbf{a}_1)\cdot \epsilon_1(\mathbf{a}_2)$, which settles the claim.
\end{proof}

\bfact
A first attempt at making this construction without the requirement of semisimplicity can be found in \cite{LeuTPMV}. The main proof was based on \cite[Theorem 4.11]{LeuTens}, which turned  out to contain a mistake, see  \cite[Remark 3.1]{LLTP1} and \Cref{sec:SEP} for more details.
\efact

\section{A category-theoretical point of view} \label{sec:CatSett}
In this section we complete the internal square of Figure 1 and finally connect all the expansions of MV-algebras involved in this work.

\subsection{From MV-algebras to PMV-algebras}
Let us recall that $\mathbf{MV_{ss}}$ is the full subcategory of MV-algebras whose objects are semisimple MV-algebras, while $\mathbf{uPMV_{ss}}$ is the full subcategory of PMV-algebras whose objects are unital and semisimple (and therefore commutative) PMV-algebras.

We define a functor $\mathbf{T}: \mathbf{MV_{ss}} \rightarrow \mathbf{uPMV_{ss}}$ as follows.
\begin{itemize}
\item[(i)] For any $A\in \mathbf{MV_{ss}}$, $\mathbf{T}(A)$ is the tensor PMV-algebra $T(A)$. By Proposition \ref{pro:tenPMV} it is a unital and semisimple PMV-algebra.
\item[(ii)] For any homomorphism of MV-algebras $h:A\rightarrow B$, $\mathbf{T}(h)$ is the homomorphism of PMV-algebras $h^{\sharp}$ defined in Corollary \ref{pro:tenhom}.
\end{itemize}
\noindent From $\mathbf{uPMV_{ss}}$ to $\mathbf{MV_{ss}}$ we have the usual forgetful functor $\mathcal{U}_{(\cdot)}$.

\blem \label{lem:TPMVfunctor}
The above-defined  map $\mathbf{T}$ is a functor.
\elem
\begin{proof}
Denoted by $\mathbf{I}_A$ and $\mathbf{I}_{\mathbf{T}(A)}$ the identity maps on $A$ and $\mathbf{T}(A)$ respectively, it is easy to check that $\mathbf{I}_{\mathbf{T}(A)} \circ \epsilon_{1,A} = \epsilon_{1,A} \circ \mathbf{I}_A$, therefore $\mathbf{I}_A^{\sharp}= \mathbf{I}_{\mathbf{T}(A)}$.\\
Let $h:A\rightarrow B$ and $g: B\rightarrow C$ be homomorphisms of  MV-algebras. We have $(g^{\sharp} \circ h^{\sharp})\circ \epsilon_{1,A} = g^{\sharp} \circ (h^{\sharp}\circ \epsilon_{1,A})=g^{\sharp} \circ (\epsilon_{1,B} \circ h)=(g^{\sharp} \circ \epsilon_{1,B}) \circ h= \epsilon_{1,C} \circ (g\circ h)$, 
then $(g^{\sharp} \circ h^{\sharp})= (g\circ h)^{\sharp}$ and $\mathbf{T}$ is a functor.
\end{proof}

\blem \label{lem:adj}
The maps $\{ \epsilon_{1,A} \} _{A\in \mathbf{MV_{ss}}}$ are a natural transformation between the identity functor on $\mathbf{MV_{ss}}$ and the composite functor $\mathcal{U}_{(\cdot)}\circ \mathbf{T}$.
\elem
\begin{proof}
Let $h: A\rightarrow B$ be a homomorphism of MV-algebras. We need to prove that $(\mathcal{U}_{(\cdot)}\circ \mathbf{T})(h) \circ \epsilon_{1,A} = \epsilon_{1,B} \circ h$. Since $(\mathcal{U}_{(\cdot)}\circ \mathbf{T})(h) =h^{\sharp}$ the result follows from Corollary \ref{pro:tenhom}.
\end{proof}

\bthm \label{teo:adjMVPMV}
The functors $\mathbf{T}$ and $\mathcal{U}_{(\cdot)}$ are an adjoint pair of functors.
\ethm
\begin{proof}
In order to prove that $\mathbf{T}$ is a left adjoint functor for $\mathcal{U}_{(\cdot)}$, we need to prove that for any unital and semisimple PMV-algebra $P$ and any homomorphism of MV-algebras $f: A\rightarrow \mathcal{U}_{(\cdot)}(P)$, with $A\in \mathbf{MV_{ss}}$, there exists a homomorphism of PMV-algebras $f^{\sharp}: \mathbf{T}(A)\rightarrow P$ such that $\mathcal{U}_{(\cdot)}(f^{\sharp})\circ \iota_A = f$. This follows from Theorem \ref{teo:tenUnProp} and Lemma \ref{lem:adj}.
\end{proof}

\subsection{From Riesz MV-algebras to \textit{f}MV-algebras}
Let us recall that $\mathbf{RMV_{ss}}$ is the full subcategory of Riesz MV-algebras whose objects are semisimple Riesz MV-algebras, while $\mathbf{ufMV_{ss}}$ is the full subcategory of \textit{f}MV-algebras whose objects are unital and semisimple (and therefore commutative) \textit{f}MV-algebras.

We define a functor $\mathcal{F_{\mathbf{T}}}: \mathbf{RMV_{ss}} \rightarrow \mathbf{ufMV_{ss}}$ as follows.
\begin{itemize}
\item[(i)] For any $R\in \mathbf{RMV_{ss}}$, $\mathcal{F_{\mathbf T}}(R)$ is the tensor PMV-algebra $T(R)$. By Theorem \ref{thm:tenfMV} it is a unital and semisimple \textit{f}MV-algebra.
\item[(ii)] For any homomorphism of Riesz MV-algebras $h:R_1\rightarrow R_2$, $\mathcal{F}(h)$ is the homomorphism $h^{\sharp}$ defined in Proposition \ref{pro:tenhom}. It is a homomorphism of \textit{f}MV-algebras by \cite[Corollary 3.11]{LeuRMV}.
\end{itemize}
\noindent From $\mathbf{ufMV_{ss}}$ to $\mathbf{RMV_{ss}}$ we have the usual forgetful functor $\mathcal{U}_{(\cdot)}$.

\bthm \label{teo:adjRMVfMV}
The two functors $\mathcal{F_{\mathbf T}}$ and $\mathcal{U_{(\cdot)}}$ constitute an adjoint pair of functors. The maps $\{ \epsilon_{1,R} \} _{R\in \mathbf{RMV_{ss}}}$ are a natural transformation between the identity functor on $\mathbf{RMV_{ss}}$ and the composite functor $\mathcal{U_{(\cdot)}}\circ \mathcal{F}$.
\ethm
\begin{proof}
It is similar to the proofs of Lemma \ref{lem:TPMVfunctor}, Lemma \ref{lem:adj} and  Theorem \ref{teo:adjMVPMV}.
\end{proof}

\subsection{Closing the circle}\label{dia}
In Section \ref{sec:intrTP} we have described the following adjunctions:

1) ($\mathcal{T_{\ten}}$, $\mathcal{U}_{\mathbb R}$) between semisimple MV-algebras and semisimple Riesz MV-algebras, 

2) ($\mathcal{F_{\ten}}$, $\mathcal{U}_{\mathbb R}$) between unital and semisimple PMV-algebras and unital and semisimple \textit{f}MV-algebras.

Together with the results of Theorem \ref{teo:adjMVPMV} and Theorem \ref{teo:adjRMVfMV}, we have two different paths from semisimple MV-algebras to unital and semisimple \textit{f}MV-algebras, as displayed in Figure 4.
\begin{center}
\begin{tikzpicture}
  \node (A) {$\mathbf{MV_{ss}}$};
  \node (B) [below of=A] {$\mathbf{uPMV_{ss}}$};
  \node (C) [right of=B] {$\mathbf{ufMV_{ss}}$};
  \node (D) [right of=A] {$\mathbf{RMV_{ss}}$};
  \draw[->] (A) to node [swap]{$\mathbf{T}$} (B);
  \draw[->] (B) to node {$\mathcal{F_{\ten}}$} (C);
  \draw[->] (A) to node {$\mathcal{T_{\ten}}$} (D);
  \draw[->] (D) to node {$\mathcal{F_{\mathbf T}}$} (C);
  \node (K)  at ($(B)!0.5!(C)$) { };
  \node [below=0.5cm, align=flush center] at (K){ Figure 4. From MV-algebras to \textit{f}MV-algebras.};

 \end{tikzpicture}
\end{center}

\bthm \label{adjunctionsCommute}
The functors $\mathcal{F_{\mathbf{T}}}\circ \mathcal{T_{\ten}}$ and $\mathcal{F_{\ten }}\circ \mathbf{T}$ are naturally isomorphic. 
\ethm

\begin{proof}
It follows by well-known properties of the adjoint functors \cite{maclane}, since  the functors $\mathcal{F_{\mathbf{T}}}\circ \mathcal{T_{\ten}}$ and $\mathcal{F_{\ten }}\circ \mathbf{T}$ are both adjoints of the forgetful functor from 
$\mathbf{ufMV_{ss}}$ to $\mathbf{MV_{ss}}$. A direct proof can be found in \cite{LPhD}.
\end{proof}

\bfact Note that we have obtained an adjunction between $\mathbf{MV_{ss}}$ and $\mathbf{ufMV_{ss}}$. If $A$ is a semisimple MV-algebra then $[0,1]\ten{\mathbf T}(A)\simeq {\mathcal{F_{\mathbf{T}}}}([0,1]\ten A)$ and this will be called the {\em semisimple tensor 
\textit{f}MV-algebra} of $A$.  
\efact

\subsection{A direct application: the amalgamation property} \label{sec:amalg}
The results of Section \ref{sec:TPMV} are useful not only in giving a complete picture of the hierarchy of MV-algebras with product, but they can be used to infer new relevant results on the classes of algebras involved. An example of this is the amalgamation property for semisimple Riesz MV-algebras, unital and semisimple PMV-algebras and unital and semisimple  \textit{f}MV-algebras.

\bprop \label{amalg:PMV}
$\mathbf{uPMV_{ss}}$, $\mathbf{ufMV_{ss}}$ and $\mathbf{RMV_{ss}}$ have the amalgamation property. 
\eprop
\begin{proof} We give the complete proof for $\mathbf{uPMV_{ss}}$.
Let $A$, $B$, $Z$  be unital and semisimple  PMV-algebras such that $Z$ embeds in both $A$ and $B$, with embeddings $z_A$ and $z_B$. We have to prove that there exists a unital and semisimple PMV-algebra $E$ such that both $A$ and $B$ embed in $E$, with embeddings $f_A$ and $f_B$ and $f_B \circ z_B=f_A \circ z_A$.

We consider the MV-reducts of $A$, $B$ and $Z$. By \cite[Theorem 2.20]{MunBook}, there exists a MV-algebra $C$ and $f_A, f_B$ embeddings such that $f_A:A \hookrightarrow C$, $f_B: B \hookrightarrow C$. It is easy to prove that, since we deal with semisimple algebras, $\pi \circ f_A$ and $\pi \circ f_B$ are embeddings of $A$ and $B$ respectively in $D=C/Rad(C)$, where $\pi:C \to C/Rad(C)$ is the canonical epimorphism.\\
By Lemma \ref{lem:emb1} $D=T^1(D)$ embeds in $T(D)$ with embedding $\epsilon_{1,D}$. We get two embeddings $\overline{f_A}:A \hookrightarrow T(D)$ and $\overline{f_B}:B \hookrightarrow T(D)$, where $\overline{f_A}=\epsilon_{1,D} \circ \pi \circ f_A$ and $\overline{f_B}=\epsilon_{1,D} \circ \pi \circ f_B$.
\begin{center}
\begin{tikzpicture}
  \node (A) {$A$};
  \node (E) [below=1cm of A] {$B$};
  \node (K) at ($(A)!0.5!(E)$) { };
  \node (L) [left of=K] {$Z$};
  \draw[right hook->] (L) to node {$z_A$} (A);
  \draw[right hook->] (L) to node [swap]{$z_B$} (E);
  \node (B) [right=1cm of K] {$D=\quot{C}{Rad(C)}$};
  \draw[right hook->] (A) to node {$\pi \circ f_A$} (B);
  \draw[right hook->] (E) to node [swap]{$\pi \circ f_B$} (B);
  \node (C) [right=2cm of B] {$T(D)$};
  \draw[right hook->] (B) to node {$\epsilon_{1,D}$} (C);
  \node (Y) at ($(B)!0.5!(C)$) { };
  \node [below=1.5cm, align=flush center] at (B){Figure 5. Amalgamation property.};
 \end{tikzpicture}
\end{center}
By \cite[Theorem 2.20]{MunBook} $f_A \circ z_A = f_B \circ z_B$, therefore the diagram commutes. Since $A$ and $B$ are unital and commutative, the corresponding \textit{f}-rings are Archimedean.  By \cite[Theorem 1.1]{BH}, any homomorphism of $\ell$-groups between unital and Archimedean \textit{f}-rings commutes with the product, and we infer that $\overline{f_A}$ and $\overline{f_B}$ are homomorphisms of PMV-algebras, settling the claim.

The proofs for $\mathbf{ufMV_{ss}}$ and $\mathbf{RMV_{ss}}$  are similar:  we shall further embed $T(D)$ in $[0,1]\ten T(D)$ in the case of \textit{f}MV-algebras and $D$ in $[0,1]\ten D$ in the case of Riesz MV-algebras. Moreover,  we shall use  \cite[Proposition 3.2]{LLfMV} and, respectively \cite[Corollary 2]{LeuRMV}, in order to get the intended type of morphisms.
\end{proof}


\section{From MV-algebras to $\ell u$-groups} \label{sec:lgroups}
We now transfer the results of the previous sections to groups and we define the \textit{tensor fu-ring} of an Archimedean $\ell u$-group, by categorical equivalence. We remark that similar adjunctions can be found in literature \--- to the best of our knowledge \--- only in the case of \emph{non-ordered} groups, rings and modules.

Let us recall some results from \cite{LLTP1}.
\bthm \label{teo:SEPgroups}
The following hold:
\begin{enumerate}[label=(\roman*)]
\item If $V$ is a Riesz space with strong unit and $G$ is an Archimedean $\ell u$-group, $V\otimes_a G$ is a Riesz space with strong unit.
\item If $R$ and $S$ are unital and Archimedean $\ell u$-rings, $R\otimes_a S$ is a unital and Archimedean $\ell u$-ring.
\item If $V$ is a unital and Archimedean \textit{fu}-algebra and $R$ is a unital and Archimedean $\ell u$-ring, $V\otimes_a R$ is a unital and Archimedean \textit{fu}-algebra.
\end{enumerate}
\ethm

Let $\Xi$ be the inverse functor of $\Gamma$  and $\Xi_{*}$ be the inverse functor of $\Gamma_{*}$, with $*\in \{ \cdot, \mathbb{R}, f\}$. For the detailed construction of $\Xi$, we refer to \cite{CDM}.

\bdfn Let $(G,u)$ an Archimedean $\ell u$-group, $A=\Gamma (G,u)$ is a semisimple MV-algebra and $T(A)$ is its tensor PMV-algebra, which is unital and semisimple.

$(R,v)=\Xi_{(\cdot)}(T(A))$ will be the unital and Archimedean  \textit{tensor \textit{fu}-ring of $(G,u)$} and will be denoted by $T(G,u)$.
\edfn
\noindent Let us remark that $\ten_a$ is associative by Theorem \ref{teo:commTP} and Corollary \ref{cor:ass01}. Hence, we can safely define $T^n(G,u)$ as $\overbrace{G \ten_a \ldots \ten_aG}^{n}$. Again by \ref{teo:commTP} we entail $$\Gamma (T^n(G,u)) \simeq T^n(A).$$
By categorical equivalence, any embedding $\epsilon_{n,m}$ extends to an embedding

$\widetilde{\epsilon}_{n,m}:T^n(G,u) \rightarrow T^m(G,u)$,\\
and any embedding $\epsilon_{n,A}$ extends to an embedding

$\epsilon_{n,G}: T^n(G,u) \rightarrow T(G,u)$, in particular, for $n=1$, $(G,u)$ embeds in $T(G,u)$.
 
\bthm \label{teo:tensalg}
The pair $(T(G,u), \epsilon_{n,G})$ is the direct limit of the direct system $(T^n(G,u), \widetilde{\epsilon}_{n,m})$. Indeed, it satisfies the following universal property:

for any \textit{f}-ring $(T,w)$ and any homomorphism $f: (G,u) \rightarrow (T,w)$ there exists $g: T(G,u)\rightarrow (T,w)$ such that $g \circ \epsilon_{1,G}=f$.
\ethm
\begin{proof}
For any \textit{f}-ring $(T,w)$ and any collection of maps $\delta_n:T^n(G,u) \rightarrow T$ such that $\delta_m \circ \widetilde{\epsilon}_{n,m}=\delta_n$   there exists a homomorphism $f:T(G,u)\to (T,w)$ that makes the diagram in Figure 6 commutative.
\begin{center}
\begin{tikzpicture}
  \node (X) {$T^n(G,u)$};
  \node (T) [right of=X] {};
  \node (Y) [right of=T] {$T^m(G,u)$ };
  \node (Z) [below=1cm of T] {$T(G,u)$};
  \node (W) [below=1cm of Z] {$(T,w)$};
  \draw[->] (X) to node {$\widetilde{\epsilon} _{n,m}$} (Y);
  \draw[->] (X) to node {$\epsilon _{n,G}$} (Z);
  \draw[->] (Y) to node [swap]{$\epsilon _{m,G}$} (Z);
  \draw[->] (X) to node [swap]{$\delta_n$} (W);
  \draw[->] (Y) to node {$\delta_m$} (W);
  \draw[->][dashed] (Z) to node [swap]{$f$} (W);
    \node [below=0.5cm, align=flush center] at (W){ Figure 6. Universal property.};
 \end{tikzpicture}
\end{center}
The existence of $f$ follows by the construction of $T(A)$ as direct limit and by categorical equivalence, therefore $(T(G,u), \epsilon_{n,G})$ is the direct limit of the direct system $(T^n(G,u), \widetilde{\epsilon}_{n,m})$.
\end{proof}

\Cref{teo:SEP,teo:commTP,teo:SEPgroups,teo:tensalg} ensure that we can safely apply $\Xi$ to the diagram from Section \ref{dia}.  Denoted by $\mathbf{auG_a}$, $\mathbf{uR_a}$, $\mathbf{uRS_a}$ and $\mathbf{fuAlg_a}$ the categories of Archimedean $\ell u$-groups, unital and Archimedean $\ell u$-rings, Archimedean Riesz spaces with strong unit and unital and Archimedean \textit{fu}-algebras respectively we get
\begin{center}
\begin{tikzpicture}
  \node (1) {$\mathbf{auG_a}$};
  \node (2) [below of=1] {$\mathbf{uR_a}$};
  \node (3) [right of=2] {$\mathbf{fuAlg_a}$};
  \node (4) [right of=1] {$\mathbf{uRS_a}$};
  \draw[->] (1) to node [swap]{$\mathbf{T_a}$} (2);
  \draw[->] (2) to node {$\mathcal{F_{\ten _\mathbf{a}}}$} (3);
  \draw[->] (1) to node {$\mathcal{T_{\ten _\mathbf{a}}}$} (4);
  \draw[->] (4) to node {$\mathcal{F_{\mathbf {T _a}}}$} (3);
  \node (K) at ($(2)!0.5!(3)$) { };
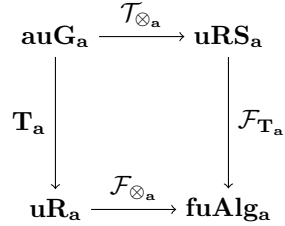
   \node [below=0.5cm, align=flush center] at (K){ Figure 7. From $\ell u$-groups to \textit{fu}-algebras.}; 
 \end{tikzpicture}
\end{center}

Moreover, applying the inverses of the functors $\Gamma$ and $\Gamma_{(\cdot)}$, $(\mathbf{T}, \mathcal{U}_{(\cdot)})$ extends to $(\mathbf{T_a}, \mathcal{U}_{\cdot \ell})$. This is an adjunction between $\mathbf{auG_a}$ and $\mathbf{uR_a}$.
Applying the inverses of the functors $\Gamma_{(\mathbb{R} )}$ and $\Gamma_{f}$ , $(\mathcal{F}_{\mathbf T}, \mathcal{U}_{(\cdot)})$ extends to $(\mathcal{F}_{\mathbf{T_a}} , \mathcal{U}_{\cdot \ell})$. This is an adjunction between $\mathbf{uRS_a}$ and $\mathbf{fuAlg_a}$.

\bfact If (G,u) is an Archimedean $\ell u$-group, we can say that  ${\mathbf {T_a}}(G,u)$ is the {\em  Archimedean tensor \textit{f}-ring} of $(G,u)$  and  ${\mathcal{F_{\ten _\mathbf{a}}}}({\mathbf {T_a}}(G,u))$ is {\em  Archimedean tensor \textit{f}-algebra} of $(G,u)$.
\efact

Finally,  Figure 8 adds the missing adjunctions to Figure 1.
\begin{center}
\begin{tikzpicture}
\node (A) {$\mathbf{auG_{a}}$};
\node (B) [below=1.5cm of A, right of=A]{$\mathbf{MV_{ss}}$};
\node (C) [right of=B] {$\mathbf{RMV_{ss}}$};
\node (D) [above=1.5cm of C, right of=C] {$\mathbf{uRS_{a}}$};
\node (E) [below of=B] {$\mathbf{uPMV_{ss}}$};
\node (F) [right of=E] {$\mathbf{ufMV_{ss}}$};
\node (G) [below=1.5cm of E, left of=E] {$\mathbf{uR_{a}}$};
\node (H) [below=1.5cm of F, right of=F] {$\mathbf{ufAlg_{a}}$};
\draw[<->] (A) to node {$\Gamma$} (B);
\draw[->] (D) to node {$\mathcal{U}_{(\ell\mathbb R)}$} (A);
\draw[->] (C) to node {$\mathcal{U}_{\mathbb{R}}$} (B);
\draw[<->] (C) to node {$\Gamma_{\mathbb R}$} (D);
\draw[->] (E) to node [swap]{$\mathcal{U}_{(\cdot)}$} (B);
\draw[->] (F) to node [swap]{$\mathcal{U}_{\mathbb R}$} (E);
\draw[->] (F) to node {$\mathcal{U}_{(\cdot)}$} (C);
\draw[->] (G) to node [swap]{$\mathcal{U}_{(\cdot \ell)}$} (A);
\draw[->] (H) to node {$\mathcal{U}_{(\cdot \ell)}$} (D);
\draw[->] (H) to node [swap]{$\mathcal{U}_{(\ell \mathbb R)}$} (G);
\draw[<->] (G) to node {$\Gamma_{(\cdot)}$} (E);
\draw[<->] (F) to node {$\Gamma_{f}$} (H);

\draw[->, bend right] (E) to node [swap]{{$\mathcal{F}_{\ten}$}} (F);
\draw[->, bend left] (B) to node {{${\mathcal T}_\ten$}} (C);
\draw[->, bend right] (B) to node [swap]{{$\mathbf{T}$}} (E);
\draw[->, bend left] (C) to node {{$\mathcal{F_{\mathbf{T}}}$}} (F);

\draw[->, bend right=12] (G) to node [swap]{{$\mathcal{F}_{\ten a}$}} (H);
\draw[->, bend left=12] (A) to node {{${\mathcal T}_{\ten a}$}} (D);
\draw[->, bend right] (A) to node [swap]{{$\mathbf{T_a}$}} (G);
\draw[->, bend left] (D) to node {{$\mathcal{F_{\mathbf{T_a}}}$}} (H);

\node (S)  at ($(G)!0.5!(H)$) { };
\node [below=1cm, align=flush center] at (S){ Figure 8. The complete diagram.};
\end{tikzpicture}
\end{center}

The following result is a straightforward consequence.

\bcor
$\mathbf{uR_a}$, $\mathbf{fuAlg_a}$  and $\mathbf{uRS_a}$ have the amalgamation property.
\ecor
\begin{proof}
It is easily deduced from Proposition \ref{amalg:PMV}  and the $\Gamma$-type categorical equivalences.
\end{proof}

\section{Description of free objects}\label{sec:free}
We close this paper providing a characterization of the free PMV-algebra, of the free \textit{f}MV-algebra and of the free Riesz MV-algebra starting from the free MV-algebra and using the tensor product. Let us point out that free algebras are related to the Lindenbaum-Tarski algebras of the corresponded logics. Thus, once again these results \--- based on algebraic and categorical methods \--- allow us to entail properties of the logical systems in a much more elegant and straightforward manner. 

We recall that the standard MV-algebra $[0,1]$ can be regarded both as a PMV-algebra and a Riesz MV-algebra, when the product and the scalar multiplication are equal with the product of  real numbers: in the sequel will denote the standard models in the obvious varieties by
\begin{align*}
&[0,1]_{MV}=([0,1], \oplus, ^*, 0)\\ 
&[0,1]_{PMV}=([0,1], \oplus, ^*, \cdot , 0)\\
&[0,1]_{fMV}=([0,1], \oplus, ^*, \cdot , \{\alpha\}_{\alpha \in [0,1]}, 0). 
\end{align*}
\bfact
Let $n\ge 1$. As stated in Section \ref{sec:intrMV}, the free PMV-algebra with $n$ free generators in
HSP$([0,1]_{PMV})$ exists and  its elements are term functions defined on $[0,1]$. More precisely, for any $n\geq 1$, assume  $X=\{x_1,\ldots, x_n\}$ and let $Term_n$ be the set of terms  with variables from $X$ in the language of PMV-algebras. We 
denote by $Free_{PMV}(n)$ the free PMV-algebra in HSP$([0,1]_{PMV})$ with $n$ free generators. It follows that
\begin{center}
$Free_{PMV}(n)=\{\widetilde{t}\colon [0,1]_{PMV}^n\to [0,1]_{PMV}\mid \widetilde{t} \text{ is the term function of } t\in Term_n\}$.
\end{center}
Since $Free_{PMV}(n)$ is a subalgebra of $[0,1]^{[0,1]^n}$ it follows that $Free_{PMV}(n)$ is unital and semisimple.
\efact
Let $Free_{MV}(n)$ and $Free_{RMV}(n)$ be the free MV-algebra and, respectively, the free Riesz MV-algebra over $n$ generators. Let $Free_{\textit{f}MV}(n)$ be the free \textit{f}MV-algebra over $n$ generators in HSP$([0,1]_{fMV})$, the variety of
\textit{f}MV-algebras generated by $[0,1]$. See more details in \cite{CDM, Mon+, LeuRMV, LLfMV}.

\bprop\label{last}
For $n\ge 1$, the following hold:
\begin{enumerate}[label=(\roman*)]
\item $Free_{RMV}(n)\simeq [0,1]_{RMV}\ten Free_{MV}(n)$,

\item  $Free_{PMV}(n)\simeq {\mathbf T}(Free_{MV}(n))$,

\item  $Free_{fMV}(n)\simeq [0,1]_{RMV}\ten {\mathbf T}(Free_{MV}(n))\simeq
{\mathcal{F_{\mathbf{T}}}}([0,1]_{RMV}\ten Free_{MV}(n))$.

\end{enumerate}
\eprop
\begin{proof}
(i) It is \cite[Proposition 5.1]{LLTP1}.

(ii)  Let $P$ be a unital and semisimple PMV-algebra and  let $f:X\to P$ be  a function, with $|X|=n$. There is a unique homomorphism of MV-algebras $f^{\#}: Free_{MV}(n)\to {\cal U}_{(\cdot)}(P)$  that extends $f$. Being free algebras semisimple, Theorem \ref{teo:tenUnProp} ensures that there exists a homomorphism of PMV-algebras $\widetilde{f}:T(Free_{MV}(n))\to P$ such that $\widetilde{f}\circ\epsilon_{1,Free_{MV}(n)}=f^{\#}$. The uniqueness of $\widetilde{f}$ is a consequence of the uniqueness of $f^\#$. Since $\epsilon_{1,Free_{MV}(n)}$ is an embedding we have 
$X\simeq \epsilon_{1,Free_{MV}(n)}(X)$ and $T(Free_{MV}(n))$ is the free object in $\mathbf{uPMV_{ss}}$. Being $Free_{PMV}(n)$ also an object in  $\mathbf{uPMV_{ss}}$, we entail that $T(Free_{MV}(n))\simeq Free_{PMV}(n)$.

(iii) It follows from (ii), \cite[Proposition 5.1]{LLTP1} and Theorem \ref{adjunctionsCommute} with similar arguments. 
\end{proof}

\paragraph*{Acknowledgement}
\noindent I. Leu\c stean was supported by a grant of the Romanian National Authority for Scientific Research and Innovation, CNCS-UEFISCDI, project number PN-II-RU-TE-2014-4-0730.

\end{document}